\documentclass[a4paper,12pt, reqno]{amsart}
\usepackage{a4wide}
\usepackage{amsthm}
\usepackage{amssymb}
\usepackage{palatino}
\usepackage{lipsum}
\usepackage{longtable}
\usepackage{enumitem}
\usepackage{hyperref}

\usepackage{rotating}
\newtheorem{thm}{Theorem}[section]

\newtheorem{prop}[thm]{Proposition}
\newtheorem{lem}[thm]{Lemma}

\theoremstyle{definition}
\newtheorem{defn}[thm]{Definition}

\theoremstyle{remark}
\newtheorem{rem}[thm]{Remark}

\makeatletter
\let\c@equation\c@thm
\makeatother
\numberwithin{equation}{section}

\begin{document}
\title[Euclidean algorithm in Galois Quartic Fields]{Euclidean algorithm in Galois Quartic Fields}
\author[K Srinivas, M Subramani and Usha K Sangale]{K Srinivas$^{(1)}$, M Subramani$^{(2)}$ and Usha K Sangale$^{(3)}$}
\address{$^{(1)}$Institute of Mathematical Sciences, HBNI,  CIT Campus, Taramani, Chennai 600 113, India}
\address{$^{(2)}$Indian Institute of Information Technology, Design and Manufacturing, Kancheepuram, Chennai 600 127, India}
\address{$^{(3)}$Institute of Mathematical Sciences, HBNI,  CIT Campus, Taramani, Chennai 600 113, India}
 \address{$^{(3)}$SRTM University, Vishnupuri, Nanded, Maharastra 431606, India}
\email{$^{(1)}$srini@imsc.res.in}
\email{$^{(2)}$subramani@iiitdm.ac.in}
\email{$^{(3)}$ushas073@gmail.com}

\begin{abstract}
We prove that all imaginary biquadratic fields and cyclic quartic fields of class number $1$ are Euclidean.
\end{abstract}

\maketitle

\section{Introduction}
Let $K$ be an algebraic number field, $\mathcal{O}_K$ be its ring of integers. We shall say that $K$ is Euclidean or $K$ admits a Euclidean algorithm, if there exists a function $\phi: \mathcal{O}_K \to \mathbb{N}\cup \{0\}$ with the property that for all $a, b\in \mathcal{O}_K, b\not= 0$ there exists  $q, r \in \mathcal{O}_K$ such that $a = bq + r$ and  $\phi(r) < \phi(b).$ If $\phi$ is the norm-function, then we say $K$ is norm-Euclidean. For example $\phi_m( r+s\sqrt{m}) = \vert r^2 - m s^2 \vert $ is a norm-function on the quadratic field $ \mathbb{Q}(\sqrt{m})$. 

\medskip

The determination of number fields which are Euclidean is an interesting topic of research in algebraic number theory. It is easy to show that a necessary condition for $K$ to be Euclidean is that it must be a principal ideal domain (PID) or equivalently they must have class number $1$. The converse is, however, not true! 
In fact, it is known that there are only  \emph{nine} imaginary quadratic fields $K = \mathbb{Q}(\sqrt{-m}),$ corresponding to $ m = 1, 2, 3, 7, 11, 19, 43, 67, 163$ with class number $1$, out of which only the first \emph{five} are Euclidean (see \cite{alaca}). On the other hand, in the case of real quadratic fields $K = \mathbb{Q}(\sqrt{m})$, it is known that $K$ is norm-Euclidean exactly when $ m = 2, 3, 5, 6, 7, 11, 13, 17, 19, 21, 29, 33,\\ 37, 41, 57, 73 $ (\cite{Chatland}, \cite{Inkeri}).  D. A. Clark \cite{CL} showed that $\mathbb{Q}(\sqrt{69})$ is Euclidean by constructing an explicit Euclidean algorithm, though it was known that it is not Euclidean with respect to the norm-function. An elegant criteria to construct Euclidean algorithm in integral domains was given by T. Motzkin \cite{Motzkin}. 
The next big step was taken in 1973 by P. J. Weinberger \cite{Weinberger} who showed that under generalized Riemann hypothesis (GRH) all algebraic number fields with infinitely many units and whose ring of integers are PIDs are in fact Euclidean!
\medskip

M. Ram Murty and D. A. Clark in \cite{Clark-Murty} introduced the concept of admissible set of primes and proved the existence of Euclidean algorithm  for a class of totally real number fields. We first recall this concept below.

\medskip

\noindent
\begin{defn} \emph{Admissible set of primes.} {Assume that $\mathcal{O}_K$ has class number one.  Let $\pi_1, \dots, \pi_s \in \mathcal{O}_K$ be distinct non-associate primes.  A set of primes $\{\pi_1, \dots , \pi_s\}$ is called an \emph{admissible} set of primes if, for all $\beta = \pi_1^{a_1} \dots \pi_s^{a_s}$ with $a_i$ non-negative integers, every co-prime residue class$\pmod{\beta}$ can be represented by a unit $\varepsilon \in \mathcal{O}_K^\times.$ In other words, the  set $\{\pi_1, \dots , \pi_s\}$ is \emph{admissible} if the canonical map $\mathcal{O}_K^\times \to \big(\mathcal{O}_K/(\pi_1^{a_1} \dots \pi_s^{a_s})\big)^\times$ is surjective.}
\end{defn}

\medskip

\noindent
Furthermore,  D. A. Clark and M. Ram Murty showed that, it is enough to check the above condition with {$a_i = 2$} for all $i$ (see page 160, \cite{Clark-Murty}). 

\medskip

\noindent
For Galois extensions with large unit ranks, the following result was proved by M. Harper and M. Ram Murty.

\medskip

\noindent
\textbf{Theorem A} [M. Harper, M. Ram Murty \cite{Harper-Murty}]
\emph{Let $K/\mathbb{Q}$ be abelian of degree $n$ with $\mathcal{O}_K$ having class number $1$, that contains a set of admissible primes with $s$ elements.  Let $r$ be the rank of the unit group.  If $r+s \geq 3,$ then $\mathcal{O}_K$ is Euclidean.}

\medskip

Recently some progress has been made to show Euclidean algorithm in real quadratic fields (see \cite{ram-srini-mani}), in certain cyclic cubic fields (see \cite{srini-mani}) and for cyclic fields of higher degree the reader may refer to \cite{Lez}, \cite{gown}.
 
\medskip

The object of this paper is to study the Euclidean algorithm in imaginary biquadratic fields and in cyclic quartic fields. 
From the works \cite{Parry-Brown}, \cite{her},  \cite{Lemmer3}, \cite{Feaver} we know that there are exactly $47$ imaginary biquadratic fields with class number $1$ and  $13$ of them are known to be norm-Euclidean (see \cite{Harper-Murty}, \cite{Inkeri}, \cite{Lemmer3}, \cite{Chatland}).
On the other hand, the class number $1$ imaginary cyclic quartic fields have been classified by B. Setzer \cite{Bennett Setzer}.
Recall that the conductor of an imaginary cyclic quartic field $K$ is defined as the smallest positive integer $f$ such that $K \subseteq \mathbb{Q}(\zeta_f)$. He proved that the class number $1$ cyclic quartic fields correspond to the conductors $f= 5,13,16,29,37,53,61$.

\medskip

In this paper we prove the following theorems.

\begin{thm}\label{thm-1}
Let $K =\mathbb{Q}(\sqrt{m},\sqrt{n})$ be an imaginary biquadratic field of class number $1$. Then $K$ is Euclidean.
\end{thm}

\begin{thm}\label{thm-2}
Every imaginary cyclic quartic field with class-number $1$ is Euclidean.
\end{thm}

Note that all the above fields $K$ have unit rank $1$, therefore if we can produce an \emph{admissible} set with $2$ elements, then by Theorem A it will follow that $K$ is Euclidean.
We shall describe below how to produce an admissible set of $2$ elements and how to find the elements in the admissible set.

\medskip

\begin{lem}\label{lem3}
Suppose that $\pi_1,\pi_2$ are unramified prime elements of residue class degree $1$, not lying above $2$. If the unit group $\mathcal{O}_K^\times$ maps onto $\left(\mathcal{O}_K/(\pi_1^2\pi_2^2)\right)^\times,$ then $\mathcal{O}_K^\times$ maps onto $\left(\mathcal{O}_K/(\pi_1^{a_1}\pi_2^{a_2})\right)^\times$ for any choice of nonnegative integers $a_1, a_2$.
\end{lem}

\noindent
\textbf{Proof.} See page 160 of \cite{Clark-Murty}.
\begin{rem}  By definition, the set $ \{ \pi_1,\pi_2 \}$ consisting of primes $\pi_1,\pi_2$ satisfying Lemma \ref{lem3} is an admissible set of primes.
\end{rem}

The following proposition gives a procedure to find the primes $\pi_1,\pi_2$ which depends on the structure of $\mathcal{O}_K^\times$. Note that for number fields listed in Table \ref{tab:long-2}, $\mathcal{O}_{K_j}^\times \cong C_g \times \mathbb{Z}$ where $g=4$ for $ j = 1, 2, \cdots , 6$; $ g = 2$ for $ j = 7, 8, \cdots, 13$ and for $ j = 19, 20, \cdots, 33$; $ g = 6$ for $ j =14, 15, 16, 17 $ and $18$. For cyclic quartic fields in Table \ref{tab:long-3}, $g=10$ for the number field with conductor $5$ and $g=2$ in the remaining cases.

\begin{prop}\label{prop-1}
Let $K$ be a number field with	$\mathcal{O}_K^\times \cong C_g \times \mathbb{Z}$. Let $\varepsilon$ be a unit of infinite order and $\eta$ be a unit of order $g$. Let $\mathfrak{p}_1$ and $\mathfrak{p}_2$ be two distinct, unramified prime ideals with inertial degree $1$ lying above $p_1$ and $p_2$ respectively and satisfying the following conditions
\begin{enumerate}
\item $\textrm{ord}_{\mathfrak{p}_1^2} (\varepsilon) = p_1(p_1-1)/g;$
\item $\textrm{gcd  } \big( p_1(p_1-1)/g,\: p_2(p_2-1)\big) = 1$;
\item $\textrm{gcd  } \big( p_1(p_1-1)/g,\: g \big) = 1$;
\item $\textrm{ord}_{\mathfrak{p}_1^2} ({\eta}) = g$;
\item $\textrm{ord}_{\mathfrak{p}_2^2} (\varepsilon) = p_2(p_2-1)$;
	\end{enumerate}
then $\mathcal{O}_K^\times$ maps onto $\big(\mathcal{O}_K/\pi_1^2 \pi_2^2\big)^\times$.	
\end{prop}	
\medskip
\noindent	
\textbf{Remark:} This is a general version of Theorem 2.1 \cite{ram-srini-mani} where the case $g=2$ was considered.

\medskip
\noindent
\textit{Proof.} Let $\beta := \varepsilon^{p_1(p_1-1)/g}.$ By (2) and (5), it follows that $\beta$ generates the group $(\mathcal{O}_K/\mathfrak{p}_2^2)^\times$.
On the other hand  by (1), we have  $$\beta \equiv 1 \pmod{\mathfrak{p}_1^2}.$$
As the image of $\eta \varepsilon$ lies in
$ (\mathcal{O}_K/\mathfrak{p}_2^2)^\times$, it follows
$$(\beta^k) (\eta \varepsilon) \equiv 1 \pmod{\mathfrak{p}_2^2}.$$ 
for some positive integer $k$ 

\medskip

\noindent
Now, let $\alpha := \eta \beta^k \varepsilon.$ 
Then by the above congruence, we have $$\alpha \equiv 1 \pmod{\mathfrak{p}_2^2}.$$ 
Also, $$\alpha \equiv \eta \varepsilon \pmod{\mathfrak{p}_1^2}.$$ 
Thus by (1) and (3),

\begin{align*}
\textrm{ord}_{\mathfrak{p}_1^2} (\eta \varepsilon) & = \frac{\textrm{ord}_{\mathfrak{p}_1^2} ({\eta}) \cdot \textrm{ord}_{\mathfrak{p}_1^2} (\varepsilon)}{\textrm{gcd } \big(\textrm{ord}_{\mathfrak{p}_1^2} ({\eta}), \textrm{ord}_{\mathfrak{p}_1^2} (\varepsilon) \big)}\\
& = \frac{g \cdot \frac{p_1(p_1-1)}{g}}{1} = p_1(p_1-1).
\end{align*}
That is, $\alpha$ generates the group $(\mathcal{O}_K/\mathfrak{p}_1^2)^\times.$ 

\medskip

\noindent
Let $(x,y) \in \big(\mathcal{O}_K/\mathfrak{p}_1^2\big)^\times \times \big(\mathcal{O}_K/\mathfrak{p}_2^2\big)^\times.$ Then there exist positive integers $e,f$ such that 
$$\alpha^e \equiv x \pmod{\mathfrak{p}_1^2} \quad  \mbox{and} \quad \beta^f \equiv y \pmod{\mathfrak{p}_2^2}.$$ 
Now the element $z := \alpha^e \beta^f$ maps onto the element $(x,y).$
Since $(x,y)$ was arbitrary, the canonical map takes  $\mathcal{O}_K^\times$ onto $\big(\mathcal{O}_K/\mathfrak{p}_1^2 \mathfrak{p}_2^2\big)^\times$ (by Chinese reminder theorem).
\qed

\section{Proof of Theorem \ref{thm-1}} 
We shall explain how to find the admissible set for one field in detail and to avoid repetition, we shall only give a list of admissible set of primes in all other cases (see Table \ref{tab:long-2}). This will complete the proof. The Sage Code for all the computations displayed in Table \ref{tab:long-2} and Table \ref{tab:long-3} are available at the URL \url{https://www.imsc.res.in/~srini/SageCodes/sagecode.txt}. 

\medskip

We shall work with the number field $K= \mathbb{Q}(\sqrt{-1},\sqrt{11}).$ In this case $\mathcal{O}_K^\times \cong C_4 \times \mathbb{Z}.$ 

\medskip
\noindent We note the following data: 
The generators for the group of units in $\mathcal{O}_K$ are $$\eta := -a, \quad \varepsilon := -\left((\frac{-1}{2}b + \frac{3}{2})a + \frac{1}{2}b - \frac{3}{2}\right)$$ where $a$ is a root of the polynomial $x^2+1$ and $b$ is a root of the polynomial $x^2-11.$ 

\medskip

We need  to choose two rational primes $p_1$ and $p_2$ satisfying the conditions of Proposition \ref{prop-1}.
A simple trial and error allows us to take $p_1 = 157$ and $p_2 = 5$. Clearly, the primes $p_1$ and $p_2$ satisfies the condition (2) and (3) of Proposition \ref{prop-1}.
Using Sage programme, we obtain the following prime ideal decomposition:
$$
p_1 \mathcal{O}_K = \left(\big(-b - \frac{3}{2}\big)a - \frac{1}{2}b\right) \left(\big(-b + \frac{3}{2}\big)a - \frac{1}{2}b\right) \left(\big(b - \frac{3}{2}\big)a - \frac{1}{2}b\right) \left(\big(b + \frac{3}{2}\big)a - \frac{1}{2}b\right),
$$
and
$$
p_2\mathcal{O}_K =  \left(\big(\frac{-1}{2}b - 1\big)a - \frac{1}{2}\right) \left(\big(\frac{1}{2}b - 1\big)a - \frac{1}{2}\right) \left(\big(\frac{-1}{2}b + 1\big)a - \frac{1}{2}\right) \left(\big(\frac{-1}{2}b - 1\big)a + \frac{1}{2}\right).
$$

\medskip

\noindent
As the primes $p_1$ and $p_2$ split completely in $\mathcal{O}_K$, we get a prime $\pi_1 := (-b - \frac{3}{2})a - \frac{1}{2}b$ lying above $p_1$ and another prime $\pi_2 := (\frac{-1}{2}b - 1)a - \frac{1}{2}$ lying above $p_2$ of inertial degree $1$.

\medskip

\noindent
Note that
$$
\varepsilon^{p_1} \equiv 14591 \not \equiv 1 \pmod{\pi_1^2}
$$
and
$$
\varepsilon^{(p_1 - 1)/4} \equiv 11776 \not \equiv 1 \pmod{\pi_1^2},
$$
but 
$$
\varepsilon^{p_1(p_1 - 1)/4} \equiv 1 \pmod{\pi_1^2}.
$$
The last congruence establishes condition (1) of Proposition \ref{prop-1}.

\medskip

\noindent
Now, from the congruence
$$
\varepsilon^{10} \equiv -1 \pmod{\pi_2^2}
$$
 the condition (5) is verified. Finally condition (4) follows from the definition of $\eta$ and a simple computation of $\textrm{ord}_{\pi_1^2} ({\eta})$ . We see that all the conditions of Proposition \ref{prop-1} are established. Thus by Lemma \ref{lem3} the set $\{(-b - \frac{3}{2})a - \frac{1}{2}b, (\frac{-1}{2}b - 1)a - \frac{1}{2}\}$ forms an admissible set. Therefore Theorem A implies the existence of Euclidean algorithm in  $K= \mathbb{Q}(\sqrt{-1},\sqrt{11})$.

\medskip
A similar computation gives the admissible set of primes for all the remaining biquadratic fields. This completes the proof of Theorem \ref{thm-1}.
\qed

\section{Proof of Theorem \ref{thm-2}}
As the cyclic quartic fields have unit rank $1$, it is enough to produce an \emph{admissible} set of primes with $2$ elements. The procedure to produce the required set of admissible primes is similar to the one described for biquadratic fields, therefore, we shall omit the details but only  list the \emph{admissible} set of primes.
\qed

\begin{center}
	\begin{longtable}{|l|l|l|}
		\caption{Imaginary biquadratic fields} \label{tab:long-2} \\
		
		\hline \multicolumn{1}{|c|}{\text{Fields}} & \multicolumn{1}{c|}{{primes}} &
		\multicolumn{1}{c|}{{admissible primes}} \\ \hline 
		\endfirsthead
		
		\multicolumn{3}{c}%
		{{\bfseries \tablename\ \thetable{} }} \\
		\hline \multicolumn{1}{|c|}{{$K$}} & \multicolumn{1}{c|}{{primes}}&
		\multicolumn{1}{c|}{{admissible primes}} \\ \hline 
		\endhead
		
		\hline \multicolumn{3}{|r|}{{Continued on next page}} \\ \hline
		\endfoot
		
		\hline 
		\endlastfoot

		$K_1 = \mathbb{Q}(\sqrt{-1},\sqrt{13})$ & $29,17$ & $((\frac{-1}{2}b - \frac{1}{2})a + \frac{1}{2}b + \frac{3}{2}), ((\frac{-1}{2}b - \frac{1}{2})a - 1)$ \\
		
		$K_2 = \mathbb{Q}(\sqrt{-1},\sqrt{19})$ & $5,73$ & $(\frac{1}{2}a + \frac{1}{2}b + 2), ((b - \frac{9}{2})a - \frac{1}{2}b + 3)$ \\
		
		$K_3 = \mathbb{Q}(\sqrt{-1},\sqrt{37})$ & $149,53$ & $(2a + \frac{1}{2}b + \frac{5}{2}), (-a + \frac{1}{2}b - \frac{7}{2})$ \\
		
		$K_4 = \mathbb{Q}(\sqrt{-1},\sqrt{43})$ & $13,17$ & $((\frac{-1}{2}b - 3)a - \frac{1}{2}), ((\frac{1}{2}b + 3)a + b + \frac{13}{2})$ \\
		
		$K_5 = \mathbb{Q}(\sqrt{-1},\sqrt{67})$ & $29,37$ & $((\frac{-2}b + \frac{33}{2})a + \frac{1}{2}b - 4), ((-b + \frac{17}{2})a + \frac{1}{2}b - 4)$ \\
		
		$K_6 = \mathbb{Q}(\sqrt{-1},\sqrt{163})$ & $53,173$ & $((-2b + \frac{51}{2})a + \frac{5}{2}b - 32), ((-b - \frac{25}{2})a + \frac{1}{2}b + 6)$ \\

		$K_7 = \mathbb{Q}(\sqrt{2},\sqrt{-11})$ & $23,31$ & $(-a - \frac{1}{2}b - \frac{1}{2}), (-a - \frac{1}{2}b + \frac{3}{2})$ \\
		
		$K_8 = \mathbb{Q}(\sqrt{-2},\sqrt{-11})$ & $3,59$ & $(a - \frac{1}{2}b + \frac{1}{2}), ((\frac{-1}{2}b - \frac{1}{2})a - \frac{1}{2}b + \frac{1}{2})$ \\
		
		$K_9 = \mathbb{Q}(\sqrt{-2},\sqrt{-7})$ & $11,43$ & $((\frac{-1}{2}b + \frac{1}{2})a - 1), ((\frac{1}{2}b + \frac{1}{2})a + 3)$ \\

		$K_{10} = \mathbb{Q}(\sqrt{-2},\sqrt{-19})$ & $11,17$ & $(a + \frac{1}{2}b + \frac{1}{2}), (2a - \frac{1}{2}b + \frac{1}{2})$ \\
		
		$K_{11} = \mathbb{Q}(\sqrt{-2},\sqrt{29})$ & $59, 83$ & $(-a + \frac{1}{2}b - \frac{5}{2}), (-a + \frac{1}{2}b + \frac{1}{2})$ \\
		
		$K_{12} = \mathbb{Q}(\sqrt{-2},\sqrt{-43})$ & $11,17$ & $((\frac{1}{2}b + \frac{5}{2})a + \frac{1}{2}b - \frac{9}{2}), (2a - \frac{1}{2}b - \frac{1}{2})$ \\
		
		$K_{13} = \mathbb{Q}(\sqrt{-2},\sqrt{-67})$ & $19,17$ & $(3a - \frac{1}{2}b + \frac{1}{2}), ((-b + 3)a + \frac{1}{2}b + \frac{23}{2})$ \\
		
		$K_{14} = \mathbb{Q}(\sqrt{-3},\sqrt{41})$ & $31,73$ & $((\frac{1}{4}b + \frac{7}{4})a + \frac{1}{4}b + \frac{3}{4}), ((\frac{-1}{4}b + \frac{9}{4})a + \frac{1}{4}b - \frac{5}{4})$ \\
		
		$K_{15} = \mathbb{Q}(\sqrt{-3},\sqrt{-43})$ & $31,79$ & $((\frac{-3}{4}b + \frac{19}{4})a + \frac{5}{4}b + \frac{35}{4}), ((\frac{1}{2}b + 2)a - \frac{1}{2}b + 5)$ \\
		
		$K_{16} = \mathbb{Q}(\sqrt{-3},\sqrt{-67})$ & $439,19$ & $((\frac{1}{4}b - \frac{33}{4})a + \frac{7}{4}b + \frac{17}{4}),((\frac{-1}{2}b + \frac{19}{2})a - \frac{2}{b} - 7)$ \\
		
		$K_{17} = \mathbb{Q}(\sqrt{-3},\sqrt{89})$ & $607,97$ & $((\frac{-3}{4}b - \frac{27}{4})a + \frac{5}{4}b + \frac{49}{4}), ((\frac{-1}{4}b - \frac{9}{4})a - \frac{1}{4}b - \frac{5}{4})$ \\
		
		$K_{18} = \mathbb{Q}(\sqrt{-3},\sqrt{-163})$ & $43,61$ & $((\frac{9}{4}b - \frac{37}{4})a - \frac{5}{4}b - \frac{199}{4}), ((\frac{-3}{4}b - \frac{7}{4})a + \frac{1}{4}b + \frac{67}{4})$ \\

		$K_{19} = \mathbb{Q}(\sqrt{-7},\sqrt{-11})$ & $23,37$ & $((\frac{-1}{4}b - \frac{1}{4})a + \frac{1}{4}b - \frac{3}{4}), ((\frac{-1}{4}b - \frac{1}{4})a + \frac{1}{4}b + \frac{5}{4})$ \\
		
		$K_{20} = \mathbb{Q}(\sqrt{-7},\sqrt{13})$ & $23,29$ & $(\frac{-1}{2}a - \frac{1}{2}b - 1), (\frac{-1}{2}a - \frac{1}{2}b + 2)$ \\
		
		$K_{21} = \mathbb{Q}(\sqrt{-7},\sqrt{-19})$ & $11,137$ & $(-a - \frac{1}{2}b + \frac{1}{2}), ((\frac{-1}{4}b - \frac{7}{4})a - \frac{5}{4}b + \frac{9}{4})$ \\
		
		$K_{22} = \mathbb{Q}(\sqrt{-7},\sqrt{-43})$ & $11,53$ & $((\frac{-1}{4}b - \frac{3}{4})a - \frac{1}{4}b + \frac{17}{4}), ((\frac{-1}{2}b - \frac{5}{2})a - b + 9)$ \\
		
		$K_{23} = \mathbb{Q}(\sqrt{-7},\sqrt{61})$ & $107,137$ & $((\frac{-1}{4}b - \frac{7}{4})a - \frac{3}{4}b - \frac{25}{4}), ((\frac{-1}{2}b + 4)a + \frac{1}{2})$ \\
		
		$K_{24} = \mathbb{Q}(\sqrt{-7},\sqrt{-163})$ & $43,179$ & $((\frac{-5}{4}b + \frac{53}{4})a - \frac{11}{4}b - \frac{169}{4}), (\frac{5}{2}a - \frac{1}{2}b - 1)$ \\
		
		$K_{25} = \mathbb{Q}(\sqrt{-11},\sqrt{17})$ & $47,59$ & $(\frac{-1}{2}a - \frac{1}{2}b - 1), ((\frac{1}{4}b + \frac{3}{4})a + \frac{1}{4}b + \frac{7}{4})$ \\
		
		$K_{26} = \mathbb{Q}(\sqrt{-11},\sqrt{-19})$ & $23,5$ & $((\frac{-1}{2}b + 2)a - \frac{3}{2}b - 7), ((\frac{1}{4}b - \frac{1}{4})a - \frac{1}{4}b - \frac{15}{4})$ \\
		
		$K_{27} = \mathbb{Q}(\sqrt{-11},\sqrt{-67})$ & $47,59$ & $(-a - \frac{1}{2}b - \frac{1}{2}), ((\frac{-1}{4}b + \frac{47}{4})a - \frac{19}{4}b - \frac{27}{4})$ \\
		
		$K_{28} = \mathbb{Q}(\sqrt{-11},\sqrt{-163})$ & $199,53$ & $((\frac{-3}{2}b - \frac{27}{2})a - \frac{7}{2}b + \frac{127}{2}), ({-2}a - \frac{1}{2}b + \frac{1}{2})$ \\
		
		$K_{29} = \mathbb{Q}(\sqrt{-19},\sqrt{-67})$ & $23,47$ & $(a + \frac{1}{2}b - \frac{1}{2}), ((\frac{-3}{4}b + \frac{47}{4})a + \frac{25}{4}b + \frac{107}{4})$ \\
		
		$K_{30} = \mathbb{Q}(\sqrt{-19},\sqrt{-163})$ & $43,47$ & $((\frac{-3}{4}b + \frac{3}{4})a - \frac{1}{4}b - \frac{167}{4}), ((\frac{83}{4}b + \frac{495}{4})a - \frac{169}{4}b + \frac{4619}{4})$ \\
		$K_{31} = \mathbb{Q}(\sqrt{-43},\sqrt{-67})$ & $23,17$ & $((\frac{-3}{4}b + \frac{181}{4})a - \frac{145}{4}b - \frac{161}{4}), ((\frac{-483}{4}b - \frac{1609}{4})a + \frac{1289}{4}b - \frac{25925}{4}) $ \\
		
		$K_{32} = \mathbb{Q}(\sqrt{-43},\sqrt{-163})$ & $47,53$ & $(a + \frac{1}{2}b - \frac{1}{2}),((\frac{-3}{4}b - \frac{37}{4})a - \frac{19}{4}b + \frac{251}{4}) $ \\ 
		
		$K_{33} = \mathbb{Q}(\sqrt{-67},\sqrt{-163})$ & $47,167$ & $(\frac{-33610220037656305}{4}b + \frac{327679543906041243}{4})a  $ \\
		&  & \hspace{1.5cm} $+ ( \frac{210083976696459637}{4}b + \frac{3512388602701730257}{4}),$ \\
		&  & \hspace{1cm}$(\tiny{(\frac{577065530491}{2}b - \frac{12863925015009}{2})a} $ \\
		&  & \hspace{1.5cm} $- (\frac{8247400771081}{2}b - \frac{60305418710075}{2})$ \\
	\end{longtable}
\end{center}

\begin{center}
	\begin{longtable}{|l|l|l||l|}
		\caption{Cyclic quartic fields.} \label{tab:long-3} \\
		
		\hline \multicolumn{1}{|c|}{{conductor}} & \multicolumn{1}{c|}{\textbf{$p_1,p_2$}} &
		\multicolumn{1}{c|}{\textbf{$\pi_1$}} &
		\multicolumn{1}{c|}{\textbf{$\pi_2$}} \\ \hline 
		\endfirsthead
		
		\multicolumn{3}{c}%
		{{\bfseries \tablename\ \thetable{} -- continued from previous page}} \\
		\hline \multicolumn{1}{|c|}{\textbf{conductor}} & \multicolumn{1}{c|}{\textbf{$p_1,p_2$}}&
		\multicolumn{1}{c|}{\textbf{$\pi_1$}} & \multicolumn{1}{c|}{\textbf{$\pi_2$}} \\ \hline 
		\endhead
		
		\hline \multicolumn{3}{|r|}{{Continued on next page}} \\ \hline
		\endfoot
		
		\hline 
		\endlastfoot
		
		$5$ & $11,31$ & $-\frac{1}{2}a^3 - \frac{1}{2}a^2 - \frac{5}{2}a - 1$ &  $-a^3 + a^2 - 3a + +1$ \\
		$13$ & $79,29$ & $(79, \frac{1}{6}a^2 - \frac{1}{2}a + \frac{65}{6})$ &  $(29, \frac{1}{6}a^2 - \frac{1}{2}a + \frac{23}{6})$ \\
		$16$ & $23,17$ & $a^2 - a + 3$ &  $a^3 + 4a - 1$ \\
		$29$ & $7,53$ & $(7, \frac{1}{10}a^2 + \frac{1}{2}a - \frac{3}{10})$ &  $(53, \frac{1}{10}a^2 - \frac{1}{2}a - \frac{183}{10})$ \\
		$37$ & $7,53$ & $(7, \frac{1}{24}a^3 + \frac{1}{24}a^2 + \frac{67}{24}a + \frac{67}{24})$ &  $(53, \frac{1}{24}a^3 - \frac{1}{24}a^2 + \frac{67}{24}a - \frac{451}{24})$ \\
		$53$ & $107,89$ & $(107, \frac{1}{14}a^2 - \frac{1}{2}a + \frac{527}{14})$ &  $(89, \frac{1}{14}a^2 - \frac{1}{2}a + \frac{401}{14})$ \\
		$61$ & $47,73$ & $(47, \frac{1}{78}a^3 - \frac{1}{78}a^2 + \frac{289}{78}a + \frac{538}{39})$ &  $(73, \frac{1}{78}a^3 - \frac{1}{78}a^2 + \frac{289}{78}a - \frac{8}{39})$ \\
	\end{longtable}
\end{center}

\noindent
\textbf{Acknowledgements.} The authors are immensely grateful to Prof. Wladyslaw Narkewicz and the anonymous referee for pointing out several inaccuracies and for suggesting changes at several places in an earlier version of this paper.
The work was initiated when the second author was a post doctoral fellow at Harishchandra Research Institute (HRI). He thanks HRI for providing excellent working environment. He further acknowledges with thanks the financial support received from SERB MATRICS Project no. MTR/2017/001006. We acknowledge with thanks the help received from SageMath programme.

\end{document}